6
\RequirePackage[l2tabu, orthodox]{nag}
\documentclass[english]{article}     

\usepackage[utf8]{inputenc} %
\usepackage[T1]{fontenc}    %
\usepackage{hyperref}       %
\usepackage{url}            %
\usepackage{booktabs}       %
\usepackage{amsfonts}       %
\usepackage{nicefrac}       %
\usepackage{microtype}      %

\usepackage{mathptmx}    %
\usepackage[english]{babel}

\usepackage{amsmath}
\usepackage{amssymb}
\usepackage{amsthm}
\usepackage{graphicx}
\usepackage{subcaption}
\usepackage{breakurl}
\usepackage{todonotes}
\usepackage{verbatim}
\usepackage{mathtools}
\usepackage{dsfont}
\mathtoolsset{showonlyrefs=true}
\usepackage[inline]{enumitem}
\usepackage{algorithm}
\usepackage{algpseudocode}
\usepackage{float}
\usepackage{siunitx}
\usepackage{icomma}
\usepackage{thmtools}
\usepackage{thm-restate}
\usepackage[a4paper, margin=1.4in]{geometry}

\newtheorem{theorem}{Theorem}
\newtheorem{proposition}{Proposition}
\newtheorem{definition}{Definition}
\newtheorem{assumption}{Assumption}
\newtheorem{remark}{Remark}

\newtheorem{corollary}{Corollary}

\newcommand{\EE}{\mathbb{E}}

\newcommand{\R}{\mathbb{R}}
\newcommand{\N}{\mathbb{N}}
\newcommand{\Z}{\mathbb{Z}}
\newcommand{\Prob}{\mathbb{P}}

\newcommand{\W}{\Omega}

\newcommand{\F}{\mathcal{F}}

\DeclareMathOperator{\argmin}{arg\min}

\DeclareMathOperator{\sign}{\text{sign }}

\title{Grid-Free Computation of Probabilistic Safety with Malliavin Calculus } %

\author{Francesco Cosentino\thanks{\scriptsize{Mathematical Institute, University of Oxford  \& The Alan Turing Institute, \url{name.surname@maths.ox.ac.uk}}}
\and Harald Oberhauser\thanks{\scriptsize{Mathematical Institute, University of Oxford  \& The Alan Turing Institute, \url{name.surname@maths.ox.ac.uk}}}
\and Alessandro Abate\thanks{\scriptsize{Dept. of Computer Science,  University of Oxford  \& The Alan Turing Institute, \url{name.surname@cs.ox.ac.uk}}}
}

\date{}

\begin{document}

\maketitle

\begin{abstract}
This work concerns continuous-time, continuous-space stochastic dynamical systems described by stochastic differential equations (SDE).  
It presents a new approach to compute probabilistic safety regions, namely sets of initial conditions of the SDE associated to trajectories that are safe with a probability larger than a given threshold.  
The approach introduces a functional that is minimised at the border of the probabilistic safety region, 
then solves an optimisation problem using techniques from Malliavin Calculus, which computes such region. 
Unlike existing results in the literature, 
the new approach allows one to compute probabilistic safety regions without gridding the state space of the SDE. 
\end{abstract}

\section{Background} 

In Control Engineering and in Formal Verification, 
a fundamental and common problem is safety analysis: this concerns identifying states of a dynamical model that are safe, namely that are associated to trajectories that do not escape (whether over finite or infinite time) a given set that is deemed to be safe~\cite{Abate2007,Abate2008,Blanchini2007}.  
Dually, one can express a reachability analysis problem by identifying states that are associated with trajectories entering a given target set.  
In the context of probabilistic models, such as stochastic differential equations (SDE), we are interested in characterising and computing the likelihood with which trajectories of the stochastic process either stay within a given set, or dually reach a target set - the former has been often studied in probability theory as the \emph{exit time} problem. Alternatively, for stochastic models we might be interested in computing the set of initial conditions associated with dynamics that are safe with a probability at least equal to, say $p$ - this is also known as $p$-safe analysis or computation of the $p$-safe region \cite{DBLP:conf/cdc/WisniewskiB17}. 

In this work, we present a new application of Malliavin Calculus \cite{Malliavin2006} to the computation of the $p$-safe region borrowing ideas from Mathematical Finance: in particular, we leverage and tailor techniques for the computation of the (so called) \emph{Greeks of a derivative} \cite{Hull2017} for our goal.  
This allows one to compute probabilistic safety regions without gridding the state space of the SDE: grid-based techniques are by-and-large the standard approach in existing literature, with known limits related to precision and computational scalability.  

\section{Related literature on Probabilistic Safety} 
Safety analysis, a standard specification in Formal Verification, has been studied on dynamical models within the Hybrid Systems community \cite{Blanchini2007}. 
Corresponding safety notions for stochastic models (and in particular for stochastic and hybrid ones - the latter feature is not under study in this work) have been explicitly introduced only over the past two decades \cite{Abate2008}, as further surveyed next. 

This work, unlike \cite{Abate2008}, focuses on continuous-time models:
particularly relevant for this setup, \cite{DBLP:conf/hybrid/HuLS00} has presented a new modeling framework named \emph{stochastic hybrid system} (SHS), 
which extends with randomness the deterministic framework of hybrid models by allowing the continuous flow inside each invariant set of the discrete state variables to be governed by stochastic differential equations (SDE), rather than deterministic ODEs. 
\cite{DBLP:conf/hybrid/HuLS00} proposes the notion of embedded Markov chain (EMC) and studies the exit probability problem, which is related to reachability analysis: 
it is shown that this probability over the EMC converges to its counterpart for the original SHS, as the discretisation step goes to zero.   
\cite{BJ06} blends the models from \cite{DBLP:conf/hybrid/HuLS00} with Markov models with jumps in \cite{d1993}, 
setting up \emph{Markov strings} and thus obtaining a very general class of models for SHS. 
Closely related to \cite{BJ06}, \cite{DBLP:conf/hybrid/Bujorianu04} introduces a general model for SHS and proposes a method based on {Dirichlet} forms, 
to study the reachability problem over SHS models. 
Similarly over SHS, \cite{DBLP:conf/hybrid/KoutsoukosR06} proposes a method to compute probabilistic reachability: 
underpinned by seminal work in \cite{ST02}, \cite{DBLP:conf/hybrid/KoutsoukosR06} first shows that reachability can be characterised as a viscosity solution of a system of coupled Hamilton-Jacobi-Bellman equations; 
second, it presents a numerical method for computing the solution based on discrete approximations, 
showing that this solution converges to the one for the original SHS model as the discretisation becomes ever finer. 
\cite{DBLP:conf/hybrid/RileyKR09} builds on \cite{DBLP:conf/hybrid/KoutsoukosR06} by employing Monte Carlo (MC) techniques for estimating probabilities of events, 
and \cite{DBLP:conf/hybrid/RileyKR09} uses multilevel splitting (MLS), a variance-reduction technique that can improve both efficiency and accuracy. 
Again over SHS, 
\cite{MOHAJERINESFAHANI201643} establishes a connection between stochastic \emph{reach-avoid} problems - problems encompassing both reachability and safety, also known as \emph{constrained reachability} problems - and optimal control problems involving discontinuous payoff functions. 
Focusing on a particular stochastic optimal control problem, namely the exit-time problem mentioned above, \cite{MOHAJERINESFAHANI201643} provides its characterisation as a solution of a partial differential equation in the sense of viscosity solutions, along with Dirichlet boundary conditions. 
\cite{DBLP:conf/cdc/WisniewskiB17} establishes an optimisation scheme for computing probabilistic safety of SHS, combining the use of barrier certificates and of potential theory. 

\cite{Wetal13} presents a method to compute \emph{protection certificates}, which are closely related to the concept of $p$-safe region, elaborated later.  
As discussed in  Remark~\ref{rem:comp_Wisniewski}, \cite{DBLP:journals/tac/WisniewskiBS20,WS15} compute the $p$-safe region based on the extended generator of stochastic dynamical systems; these contributions  characterize the safety problem as an optimization problem on the space of positive measures and then solve it via a moment-based method. %
\cite{DBLP:conf/cdc/BujorianuW19} characterizes the p-safe regions using concepts from Potential Theory. %

Alternative techniques leveraging randomised approaches have been presented: \cite{DBLP:conf/eucc/HuP03} introduces a method for estimating the probability of conflict for two-aircraft encounters at a fixed altitude - a probabilistic safety problem.  
The procedure is based on the introduction of a Markov chain approximation of the stochastic process describing the relative position of the aircraft. 
Along similar lines, \cite{DBLP:journals/tits/PrandiniHLS00} discusses the maximum instantaneous probability of conflict: randomised algorithms are introduced to efficiently estimate this measure of criticality and to provide quantitative bounds on the level of the approximation introduced. Also, approximate closed-form analytical expressions for the probability of conflict are obtained. 
These randomised approaches can be related to statistical model checking (SMC) techniques, 
which have also been developed for models related to SHS in \cite{SZ16}. 

Finally, the work in \cite{DBLP:journals/tac/ZamaniEMAL14} enables sound verification and correct-by-construction controller synthesis for stochastic models and their hybrid extensions \cite{DBLP:conf/eucc/ZamaniEAL13}: a stochastic control model satisfying a probabilistic variant of \emph{incremental input-to-state stability} is shown to be abstracted into a finite-state transition system, which is epsilon-approximately bisimilar to the original model.

\section{Problem Statement}

Let us consider a $d$-dimensional Brownian motion $W_t \in \R^d$ defined on a filtered probability space $(\W, \F, (\F_t), \Prob)$, and the following SDE
\begin{align}\label{eq:SDE_ch6sec2}
dX_{t} = & \mu(X_{t})dt+\sum_{k=1}^{d}\sigma_k(X_{t})dW_{t}^{k}, \quad\quad\quad X_0=x.
\end{align} 
The setup above is 
adopted by related literature, as surveyed above. 

We introduce the following requirements, 
which are used in \cite{Malliavin2006} and in particular are sufficient for all the results and algorithms proposed in this work. Obtaining weaker requirements, and thus generalising our setup, would require modifying the technical results from Malliavin calculus, which is not core to our contributions. 
\begin{assumption}\label{ass:greeks}
We suppose that the vector fields $\mu, \sigma_k \in C^\infty_{l, b}(\R^d;\R^d), k=1, 2, \ldots, d$, where $C^\infty_{l, b}(\R^d;\R^d)$ indicates the space of infinitely differentiable functions with bounded derivatives and bounded linear growth from $\R^d$ to $\R^d$. 
Moreover, if we call $\sigma\in C^\infty_{l, b}(\R^d;\R^{d\times d})$ the matrix whose columns are the vectors $\sigma_k, k\geq1$, we assume that $\sigma$ satisfies the uniform ellipticity condition, i.e. $\sigma\sigma^\top$ is uniformly positive definite.
\end{assumption}
If Assumption~\ref{ass:greeks} holds, then it is well-known that the SDE~\eqref{eq:SDE_ch6sec2} 
has a unique strong solution $X^x_t$~\cite{NIkeda2014}, and whenever clear from the context we shall omit the index $x$.

Let us consider a bounded and %
smooth region $A$ and let $\partial A$ denote  the border of $A$. We call $\tau_A^{x}$ the exit time of $X_t^x$ from the region $A$, i.e. 
\[
\tau_{A}^{x} = \inf\{t\geq 0\,:\,X_{t}^{x}\not\in A\}, \quad\quad x\in A.
\]
Whenever clear from the context we shall omit the indexes $x, A$. 

We define $\mathbb{A}^p_T$ to be the \emph{$p$-safe portion of a region} $A$, or equivalently \emph{$p$-safe region of }$A$, as the initial points $x$ in $A$ such that if $X_t$ starts from $x$, then it stays in $A$ longer than $T$ with probability greater than $p$, i.e.%
\[
\mathbb{A}^p_T = \{x\in A : \Prob(\tau^x\geq T)\geq p\}.
\]
Again, whenever clear from the context we shall omit the indexes $p, T$. 
\begin{remark}\label{rem:comp_Wisniewski}

In \cite{DBLP:journals/tac/WisniewskiBS20} the authors study a more general problem, namely the probabilistic \textit{reach-avoid} problem, defined next.  
Given a safe set $S$ and an unsafe set $U$, they compute the probability to leave $S$ before entering in $U$, before a pre-specified time $T$, i.e. 
$\{x\in S: \Prob(\tau^x_{U^C}\leq\tau^x_S, \tau^x_{U^C}<T)\leq q\}$\footnote{
In \cite{DBLP:journals/tac/WisniewskiBS20} the authors compute $\Prob(\tau^x_{U^C}<\tau^x_S, \tau^x_{U^C}<T)\leq q$, whereas here we use $\leq$, as it does not change the outcome, whilst greatly simplifying the comparison  between \cite{DBLP:journals/tac/WisniewskiBS20} and this work. }.
If we consider a set $U$ s.t. $U^C= S$, then the wanted quantity becomes  $\{x\in S: \Prob(\tau^x_{S}<T)\leq q\}$, which is exactly the dual of $\mathbb{A}^p_T $ for 
$p =1-q$. An analysis of the approximation error is not presented and, since the approach is radically different from the one presented in this paper (cf. discussion in Related Literature and in the next section), a quantitative comparison between the two approaches is questionable.
\end{remark}

The standard way to compute $\mathbb{A}^p_T$ is to discretise the region $A$ and to compute the value $\Prob(\tau^x\geq T)$ at any point in the introduced grid  %
(cf. Related Literature). 
In the following instead, using ideas from Mathematical Finance and results from the Malliavin Calculus, we show how to compute $\mathbb{A}^p_T$ with a grid-free technique.  
The approach hinges on the observation that the border of $\mathbb{A}^p_T$ can be expressed as\footnote{We employ here for simplicity a quadratic function $f(x) = 0.5 ( x - p )^2$, however any other differentiable function $f$ minimised in $x=p$ is also appropriate for the task.}
\begin{align}\label{eq:optimization}
    \partial \mathbb{A}^p_T=\argmin_{x}\frac{1}{2}\left(\Prob\left(\tau_{}^{x}\geq T\right)-p\right)^{2}.
\end{align}
The main idea of this approach is thus to solve such optimisation problem:  
indeed, assuming differentiability and excluding convexity issues, we know that, setting up the recursion 
\begin{align}\label{eq:GD_step}
x_{j+1}=&x_{j}-\lambda\left(\Prob\left(\tau^{x}\geq T\right)-p\right)D_{x}\Prob\left(\tau^{x}\geq T\right),
\end{align}
then 
 $x_j \to x^\star$, where $x^\star \in\partial \mathbb{A}^p_T$, 
for $\lambda>0$ small enough.  Equation~\eqref{eq:GD_step} represents a standard Gradient Descent (GD) step. 
We remark that in principle other optimization algorithms can be used to solve the problem in~\eqref{eq:optimization}; in this work we focus on first-order gradient-based optimization procedures, of which GD is an exemplar. 
As an alternative instance to standard GD, in the case study we employ ADAM \cite{Kingma2014}, a state-of-the-art optimization procedure.   %

Using the GD in Equation \eqref{eq:GD_step} not only allows identifying the set $\mathbb{A}^p_T$: in Section \ref{sec:walk} we also provide a procedure to explore its border. Furthermore, in Theorem~\ref{th:noholes} we prove that the ``\textit{interior}'' of the obtained region delimited by $\partial \mathbb{A}^p_T$ is in $\mathbb{A}^p_T$, which implies that there is no need to check these internal points. 
Moreover, in Proposition~\ref{prop:ins_or_out} and Corollary~\ref{cor:ins_or_out} we show how to check if a point $x$ is ``\textit{inside}'' $\mathbb{A}^p_T$ without computing $\Prob\left(\tau^{x}\geq T\right)$, but only using the gradient at a specific point on $\partial \mathbb{A}^p_T$, which is generated by the optimization procedure.

\section{Grid-based vs -free approaches} 

Evidently, the GD step in Equation \eqref{eq:GD_step} depends on the two quantities $\Prob(\tau^x\geq T)$ and $D_x\Prob(\tau^x\geq T)$: it should be clear that if we can compute them (or approximations thereof) with a grid-free method, then the overall procedure will result in a grid-free computation of the $p$-safe portion of the region. An important advantage of such grid-free approach is that if the $p$-safe region $\mathbb{A}^p_T$ is not empty - however small, even if it was a zero-measure set - then it will be found. Instead, grid-based approaches (broadly all those presented in the previous section on related work) will find the set $\mathbb{A}^p_T$ only if it intersects with the introduced grid. 
As an extreme instance, if the $p$-safe region consists of only one point, %
the procedure introduced here shall find it, up to a numerical precision related to the approximation of $\Prob(\tau^x\geq T)$ and $D_x\Prob(\tau^x\geq T)$; on the contrary, this might not be possible for grid-based approaches, unless the grid is selected to intersect such point (which is usually not known beforehand) - and this is a limit holding regardless of their numerical implementation. {Another extreme case can be identified when the $p$-safe region of interest is not bounded: in such case the approach underpinning grid-based methods can be quite inefficient, whilst the grid-free based approach presented here shall converge to its border $\partial\mathbb{A}^p_T$, and explore as much of it as computationally feasible}. 

In general, a formal comparison between the two approaches can be problematic: whilst grid-free strategies search for solutions within an uncountable infinite set, grid-based procedures search for solutions within a pre-defined, possibly finite set. 

Still, we can comment on the computational complexity related to the two different approaches: suppose that we are working with a model of dimension 2, that set $\mathbb{A}^p_T$ consists of only one bounded connected region, and consider a grid over  $\gamma\Z^2$, where $\gamma$ is a scaling parameter. 
In the following, we will show that the procedure presented in this work requires a step-exploration parameter (again called $\gamma$) that can be  related to the $\gamma$ parameter of the grid: they both indicate how precise we want the approximation to be, 
see Figure~\ref{fig:directions} and Algorithm~\ref{algo:walk_border_d2}. 
The points explored by the two methods can be  quantified as $O(\operatorname{Area}(A)\gamma^{-2})$ for the grid discretisation - this is the number of points of $\gamma\Z^2$ in the overall (larger) region $A$ - and $O(\operatorname{Len}(\partial \mathbb{A}^p_T)\gamma^{-1}) + C $ for the method here presented - the first term represents the number of points we will explore on $\partial \mathbb{A}^p_T$, whilst $C$ depends on how many points we explore to arrive at the border $\partial \mathbb{A}^p_T$ from the starting point $x_\star$ we choose for the GD procedure. 
For a model in dimension $d$, we would have instead $O(\operatorname{Area}(A)\gamma^{-d})$ and $O(\operatorname{Area}(\partial \mathbb{A}^p_T)\gamma^{-d+1}) + C. $\footnote{This  is a slight abuse of notation, indeed $\operatorname{Area}(\partial \mathbb{A}^p_T)$ now represents the Lebesgue measure in $\R^{d-1}$, whist $\operatorname{Area}(A)$ represents the Lebesgue measure in $\R^{d}$.} 
As an estimate, we can see that, as $\gamma \to 0$, the order of points explored is much less with the method presented here.

\section{Malliavin Calculus for stopping times}

In Pricing Theory, a branch of Mathematical Finance, 
a classical problem is to evaluate the variation of the price of a derivative, 
in response to a change of the underlying asset price or volatility {\cite{Hull2017}}. 
These quantities are known as \emph{Greeks} and  play a core role in hedging theory. 
More precisely, given an underlying asset, whose price $P_t$ is the solution of an SDE starting in $p_0$, the price $D$ of a derivative is given as the expectation of a functional of $P_t$, i.e.
\[
D_{p_0} = \EE [f(P_T^{p_0})], \quad\quad T>0.
\]
If we call $p_0$ the initial price of the underlying asset, the \textit{Greek} representing the sensitivity with respect to the initial price is called $\Delta$, and is formally defined as 
\[
\Delta={\frac{\partial D_{p_{0}}}{\partial p_{0}}}=\frac{\partial\EE[f(P_{T}^{p_{0}})]}{\partial p_{0}}.
\]
Through Malliavin Calculus it is possible to provide explicit formulae for the \textit{Greeks}~\cite{Fournie1999, Fournie2001, Malliavin2006, Privault2003, Bally2005, Gobet2005}.  
We refer to~\cite{Gobet2003} for a computational perspective on these methods.

In our problem setup, we are interested to compute the quantity $D_{x}\Prob\left(\tau^{x}\geq T\right)$ used in \eqref{eq:GD_step}, where $\tau^{x}$ is a specific exit time related to the probabilistic safety property: 
this is challenging because it involves the derivative of a non-smooth indicator functional of the exit time.  
We should otherwise estimate this quantity  numerically, with associated unavoidable imprecision. 
Under Assumption~\ref{ass:greeks}, it is possible to show that $X^x_t$ is a.s. differentiable with respect to the starting point $x$~\cite{Kunita1984}. 
Following the notations in~\cite{Malliavin2006}, let us introduce 
$J_t = D_x X_t $ and $\boldsymbol{\mu}=D_x[\mu(x)]$, $\boldsymbol{\sigma}_k = D_x[\sigma_k(x)]\in C^\infty_b(\R^d; \R^{d\times d})$; then $J_t$ solves 
\begin{align}
dJ_{t}=&\boldsymbol{\mu}(X_{t})J_{t}dt+\sum_{k=1}^{d}\boldsymbol{\sigma}_k(X_{t})J_t dW_{t}^{k}, \quad\quad J_{0}=I_{d}. 
\end{align}
The main results we leverage is the following.
\begin{theorem}\label{th:Delta_Malliavin}\cite[Theorem 2.18]{Malliavin2006} %
If Assumption~\ref{ass:greeks}
holds true, calling $\tau^1$ the time when $\int_0^{\tau^1}\text{dist}(X_{t}, \partial A)^{-2}dt=1$, then
\begin{align}
\frac{\partial}{\partial\varepsilon}\Prob\left[\tau^{x+\varepsilon\varsigma}\geq T\right]\Big|_{\varepsilon=0}&=\EE\left[\mathds{1}_{\{\tau^{x}\geq T\}}H_{\varsigma, T}\right], \\ H_{\varsigma, T}&=\sum_{k=1}^{d}\int_{0}^{T}\beta^{k}_t\frac{\mathds{1}_{\{t<\tau^{1}\}}}{\text{dist}(X_{t}, \partial A)^{2}}dW_{t}^{k},
\end{align}
where ${\beta^k_t}\in \R^d$ is the stochastic process satisfying  
$\sum_{k=1}^{d}\beta^{k}_t\sigma_k(X_{t}) = {J_{t}\cdot\varsigma}$. %
\end{theorem}
Selecting the directions $\varsigma = e_i, i=1, \ldots, d$, 
we can obtain the gradient $D_{x}\Prob\left(\tau_{}^{x}\geq T\right)$, 
which lies at the core of our procedure in Equation \eqref{eq:GD_step}: without this result, this derivative should be estimated alternatively, for instance numerically. 
Therefore, we have that $D_{x}\Prob\left(\tau_{}^{x}\geq T\right)=\EE\left[\mathds{1}_{\{\tau^{x}\geq T\}}H_T\right]$, where 
\begin{align}\label{eq:beta_H}
H_T &=\int_{0}^{T}\frac{\mathds{1}_{\{t<\tau^{1}\}}}{\text{dist}(X_{t}, \partial A)^{2}}\beta_t\cdot dW_{t}, \\
\beta_t &=\sigma^{-1}(X_{t})\cdot J_{t},
\end{align}
and $\sigma$ is the matrix whose columns are the vectors $\sigma_{k}$. Please note that the dimension of $\beta_t, H_t$ in Equation~\eqref{eq:beta_H} and Theorem~\ref{th:Delta_Malliavin} are different. 

\section{Properties of the region $\mathbb{A}^p_T$}

Whilst Theorem~\ref{th:Delta_Malliavin} can be useful for the problem at hand, from an algorithmic point of view there are still a few subtle points to be handled. 

Firstly, we do not know whether the quantity $\Prob\left(\tau_{}^{x}\geq T\right)$ is convex or not. 
Nevertheless, we know that $x$ is in $\partial\mathbb{A}^p_T$ if 
$
\Prob\left(\tau^{x}\geq T\right)=p, 
$
which implies that the quantity in \eqref{eq:GD_step} 
$\left(\Prob\left(\tau^{x}\geq T\right)-p\right)D_{x}\Prob\left(\tau^{x}\geq T\right)=0$, 
regardless of the value of the gradient $D_{x}\Prob\left(\tau^{x}\geq T\right)$. Therefore, if we end at a point $x$ where
\[
\Prob\left(\tau^{x}\geq T\right)\not=p\text{ and }\left(\Prob\left(\tau^{x}\geq T\right)-p\right)D_{x}\Prob\left(\tau^{x}\geq T\right)=0, 
\]
then we know that $x\not\in\partial\mathbb{A}^p_T$, thus we are in a local saddle or local maximum point. 

Secondly, we observe that the GD scheme in \eqref{eq:GD_step} converges to a point, 
however in general it does not ``discover'' the entire border $\partial \mathbb{A}^p_T$.  
Besides,  if $\mathbb{A}^p_T$ is the union of two (or more) disconnected regions, then the GD scheme will converge solely to one of them. 
The former issue can be mitigated algorithmically, by finding a way to ``explore'' the border defined by the condition $\{\Prob\left(\tau^{x}\geq T\right)=p\}$: 
this is discussed in the next Section.
However, we cannot in general solve the latter problem, which is related to the issue of convergence to local-vs-global optima, which is intrinsic to GD schemes. 

Still, we shall shed some further light on the shape of $\mathbb{A}^p_T$. 
Let us start noticing that if $x\in A$, then $\mathbb{P}\left(\tau^{x}\geq T\right)\geq 0 $, therefore for any $p>0, T>0$ 
s.t. $\mathbb{A}_{T}^{p}\not=\emptyset$,
\[
\mathbb{A}_{T}^{p} 
\subseteq 
A=\mathbb{A}_{T}^{0}. 
\] 
However, we cannot be sure that the $p$-safe region $\mathbb{A}^p_T$ is a connected set, as we can in general  argue that $\mathbb{A}^p_T=\cup_i \mathbb{A}_i$, namely $\mathbb{A}^p_T$ consists possibly of a countably infinite union of sets, wherein any $\mathbb{A}_i$ is a bounded connected set. Each component $\mathbb{A}_i$ is endowed with interesting properties.  

\begin{definition}
We say that a surface (see \cite{Kinsey1993} for a formal definition) 
is closed if it partitions the space, e.g. $\R^d$, into one bounded connected region and one unbounded region.
We denote this bounded region as the \emph{interior} of the surface. 
\end{definition}

\begin{theorem}[No holes]\label{th:noholes}
Let the $\partial \mathbb{A}_i$ be a closed surface such that 
$\partial \mathbb{A}_i\subseteq \partial \mathbb{A}^p_T$. Then, the interior of $\partial \mathbb{A}_i$ is in $\mathbb{A}^p_T$. 
\end{theorem}
\begin{proof}
Let us indicate with $\mathbb{A}_i$ the interior of $\partial \mathbb{A}_i$. 
We prove the thesis if %
for any $x\in \mathbb{A}_i$, $\Prob(\tau_A^x\geq T)\geq p$ -- we omit the index $A$ in the next steps. 
If we define $\theta$ to be the exit time from $\mathbb{A}_i$, then
\begin{align}
\mathbb{P}\left(\tau^{x}\geq T\right)=&\mathbb{P}\left(\tau^{X_{\theta}^{x}}\geq T-\theta^{x}\Big|\theta^{x}\leq T\right)\mathbb{P}\left(\theta^{x}\leq T\right)+\mathbb{P}\left(\tau^{x}\geq T\Big|\theta^{x}\geq T\right)\mathbb{P}\left(\theta^{x}\geq T\right)\\=&\mathbb{P}\left(\tau^{X_{\theta}^{x}}\geq T-\theta^{x}\Big|\theta^{x}\leq T\right)\mathbb{P}\left(\theta^{x}\leq T\right)+\mathbb{P}\left(\theta^{x}\geq T\right), 
\end{align}
where thanks to the definition of $\theta^x, \tau^x$ we have that $\mathbb{P}\left(\tau^{x}\geq T \,| \, \theta^{x}\geq T\right)=1$, indeed $\theta^x\leq\tau^x$ a.s. since $ \mathbb{A}_i \subseteq A$ and by definition of exit time. \\
Since $\theta^x\ge 0 $ a.s., 
$\mathbb{P}(\tau^{X_{\theta}^{x}}\geq T-\theta^{x}\,|\, \theta^{x}\leq T)$ $\geq$ $\mathbb{P}(\tau^{X_{\theta}^{x}}$ $\geq T\,|\, \theta^{x}\leq T)
$ then 
\begin{align}
\mathbb{P}\left(\tau^{x}\geq T\right)\ge&\mathbb{P}\left(\tau^{X_{\theta}^{x}}\geq T\Big|\theta^{x}\leq T\right)\mathbb{P}\left(\theta^{x}\leq T\right)+\mathbb{P}\left(\theta^{x}\geq T\right)\\\geq&p\mathbb{P}\left(\theta^{x}\leq T\right)+p\mathbb{P}\left(\theta^{x}\geq T\right)\\=&p,
\end{align}
because $X_{\theta}^{x}\in\partial\mathbb{A}^p_T.$ \hfill \qedsymbol
\end{proof}

\smallskip

From an algorithmic point of view, Theorem~\ref{th:noholes} is remarkable: once the algorithm has obtained a \emph{closed} surface for $\{x : \mathbb{P}\left(\tau^{x}\geq t\right) = p\}$ we know that all the points inside are in $
\mathbb{A}^p_t$ without the need to check any further. 
Nevertheless, let us recall that we cannot know if this is the only part of $\mathbb{A}^p_T$ as there could be other bounded sets in $A$, not connected with the one just found. 

\smallskip

Once we have identified (part of) $\mathbb{A}^p_T$, an important question is how to check if a point lies inside $\mathbb{A}^p_T$.
There are different ways to check if a point is inside a region, %
such as the winding number, or the Point-in-Polygon algorithm~\cite{Haines1994,Hormann2001,Kumar2018,Huang1997}, but computationally these methods are quite expensive and generalizations to dimensions greater than $3$ do not seem to be treated in the literature, at least from an algorithmic point of view.

Remember that to compute $\partial\mathbb{A}^p_T$, we use a gradient-based optimization algorithm, requiring the computation of the quantity  $D_{x}\Prob(\tau^{x}\geq T)$ for any point in the sequence \eqref{eq:GD_step}.  
Hence, it would be useful to understand if one point is inside the safety region using the information given by $D_{x}\Prob(\tau^{x}\geq T)$:   this is handled by the next result.

\begin{proposition}\label{prop:ins_or_out}
Let us suppose that a region ${A}\in \R^d$ is defined by a differentiable function $\alpha$, i.e. ${A}:=\{x : \alpha(x)\leq 0 \}$ and $\partial {A} := \{x : \alpha(x)= 0\}$. Moreover, let us suppose that ${A}$ is connected. 
Then, a point {$x$} is inside ${A}$ if
\[
x=x^{\star}-\left\Vert x-x^{\star}\right\Vert \frac{D_{x}\alpha(x)\big|_{x=x^{\star}}}{\left\Vert D_{x}\alpha(x)\big|_{x=x^{\star}}\right\Vert }, 
\]
where $x^\star := \argmin_{y\in \partial {A}} \|x-y\|$. If instead 
\[
x=x^{\star}+\left\Vert x-x^{\star}\right\Vert \frac{D_{x}\alpha(x)\big|_{x=x^{\star}}}{\left\Vert D_{x}\alpha(x)\big|_{x=x^{\star}}\right\Vert },
\]
then $x$ is outside.
\end{proposition}
\begin{proof}
Let us consider $S = S(x, \|x-x^\star\|)$ the open sphere with center $x$ and radius $\|x-x^\star\|$; we know that if $x$ is in ${A}$ then $S\subset {A}$, vice versa $S\subset {A}^C$ if $x$ is outside $A$.\\
Note that $x-x^\star$ is perpendicular to the tangential plane to $\alpha$ in $x^\star$, as it is also the gradient $D_{x}\alpha(x)\big|_{x=x^{\star}}$, therefore 
\[
x=x^{\star}\pm\left\Vert x-x^{\star}\right\Vert \frac{D_{x}\alpha(x)\big|_{x=x^{\star}}}{\left\Vert D_{x}\alpha(x)\big|_{x=x^{\star}}\right\Vert }.
\]
Since ${A}$ is connected,  $\sign\{\alpha(x)\}$ is the same for any $x\in S$ and given that $\alpha(x^\star)=0$ the sign can be deduced by the direction of the gradient, which means that if the $D_{x}\alpha(x)\big|_{x=x^{\star}}$ points to $x$ than $\alpha(x) \geq 0$ and $x\not\in {A}$, if the $-D_{x}\alpha(x)\big|_{x=x^{\star}}$ points to $x$ than $\alpha(x) \geq 0$ and $x\in {A}$. \hfill \qedsymbol
\end{proof}

\smallskip

Since we know from Theorem~\ref{th:noholes} that any portion $\mathbb{A}_i$ of $\mathbb{A}^p_T$ is connected, 
once we have found a closed surface bordering $\mathbb{A}_i$, then thanks to Proposition~\ref{prop:ins_or_out} we know how to check if a point $x$ is inside $\mathbb{A}_i$ by estimating the gradient in  $\partial\mathbb{A}_i$, which we compute during the optimization procedure. This means that we do not have to compute $\Prob(\tau^x\geq T)$. 
Unfortunately we cannot know a-priori if it is outside because we do not know beforehand whether  $\mathbb{A}^p_T$ is connected or not.  
\begin{corollary}\label{cor:ins_or_out}
Let the $\partial \mathbb{A}_i$ be a closed surface such that 
$\partial \mathbb{A}_i\subseteq \partial \mathbb{A}^p_T$ and $\mathbb{A}_i$ its interior. 
Denoting by $x^\star := \argmin_{y\in \partial \mathbb{A}_i} \|x-y\|$, then a point $x$ is inside $\mathbb{A}_i$ if
\[
x=x^{\star}-\left\Vert x-x^{\star}\right\Vert \frac{D_{x}\Prob(\tau^x\geq T)\big|_{x=x^{\star}}}{\left\Vert D_{x}\Prob(\tau^x\geq T)\big|_{x=x^{\star}}\right\Vert }.
\]
\end{corollary}
\begin{proof}
The proof follows closely Proposition~\ref{prop:ins_or_out} considering $\mathbb{A}_i$ in place of $A$. The difference is that the sign of the points in $S$ is the same for the points inside, whilst we cannot say the same if $x$ is outside the region $\mathbb{A}_i$. It could be that, if $x$ is outside $\mathbb{A}_i$, but $S\cap \mathbb{A}_j\not=\emptyset$ for some $j\not=i$. \hfill \qedsymbol 
\end{proof}

\section{Exploration of the border $\partial\mathbb{A}^p_T$}\label{sec:walk} 
In this section we provide an algorithmic solution to the problem discussed above, namely how to explore points on the border characterised by 
\[
\argmin_{x}\frac{1}{2}\left(\Prob\left(\tau_{}^{x}\geq t\right)-p\right)^{2}. 
\]
In the following we treat in detail the case of models with dimension $d=2, 3$; if the model dimension is greater than $3$, we show that the procedure can be iterated. 

\begin{algorithm}\caption{Exploration of the border $\partial\mathbb{A}^p_T$}\label{algo:walk_border_d2}
\begin{algorithmic}[1]
\State{Initialize $x_0$}
\State{Run the GD from $x_0$ up to a point $x_\star$ such that $\mathbb{P}\left(\tau^{x_\star}\geq T\right)\approx p$}
\State{$x \leftarrow x_\star$}
\State{step $\leftarrow 1$}
\While{\textbf{not} ($x\approx x_\star$ \textbf{and} step$>$step\_min)}\label{step:terminatewhile}
\State{Append $x$ to $\partial\mathbb{A}^p_T$}
\State{Move in a direction \textit{dir} perpendicular to $ D_{x}\mathbb{P}\left(\tau^{x}\geq t\right)$ }
\State{$x \leftarrow x + \gamma \times \frac{dir}{\|dir\|}$} \label{step:dir}
\State{step $\leftarrow $ step $+1$}
\EndWhile
\end{algorithmic}
\end{algorithm}

\subsection*{Dimension 2}  Algorithm~\ref{algo:walk_border_d2} explains how to proceed if $d=2$. Let us discuss the main steps of Algorithm~\ref{algo:walk_border_d2}:
\begin{enumerate}[label=\textit{(\roman*)}]
\item Line~\ref{step:terminatewhile}. If we move along the border of the region $\mathbb{A}_i$ just found, in a, say, clockwise manner, thanks to Theorem~\ref{th:noholes} we know that we can stop when we have found a closed point ($x\approx x_\star$) and all the points inside are in $\mathbb{A}_i$. 
Given that the set could be the union of different disjoint sets, we should still explore the rest of the region $A$, i.e. $A\setminus \mathbb{A}_i$.
\item Line~\ref{step:dir}. Defining $x_0 = x + \gamma \times dir$ we know that 
\begin{align}
\mathbb{P}\left(\tau^{x_{0}}\geq t\right)\approx&\,\mathbb{P}\left(\tau^{x}\geq t\right)+D_{x}\mathbb{P}\left(\tau^{x}\geq t\right)(\gamma\times dir)+\|\gamma\times dir\|^{2}\text{ERR}\\\approx&\,p+\|\gamma\times dir\|^{2}\text{ERR},
\end{align}
which means that for small $\gamma$ we do not go far from the border $\partial \mathbb{A}^p_t$. This seems the best we can do without computing further derivatives (other than the gradient). 
ERR represents the error term of a Taylor expansion. 
It is important to remark that $\gamma$ and $\lambda$ are two different parameters, which can be chosen independently, however for more insights see Section~\ref{sec:exp}.
\end{enumerate}

\subsection*{Dimension 3} In dimension 3 we can  explore the desired border along its ``sections''.   
Without loss of generality, let us suppose that the region $A$ is the sphere of center $0$ and radius $1$. 
Let us fix the discretisation parameter $\delta>0$, which is related to the error we can tolerate. 
We can discretise the first directions $x_1$ to create the planes $x_1=\pm i\delta, i\in \N, 0\leq i\leq \delta^{-1}$. 
The sections of the border are therefore the curves resulting from the intersections between the border and the considered planes.  
We thus run Algorithm~\ref{algo:walk_border_d2} constrained on any given plane $x_1=\pm i\delta$ that we are considering, see Figure~\ref{fig:experiments_full}. 
Then, the same must be done for the other directions $x_2$ and $x_3$. 
Note that, if we have already computed a ``section'', e.g. for the plane $x_1 = 0$, then this information can be very useful for the computation of the close sections, e.g. $x_1 = \pm\delta$. \\
There are two generalizations to this procedure. 
Firstly, we can consider alternative directions: instead of selecting directions $x_1, x_2, x_3$ corresponding to the vectors in the canonical basis $(e_1, e_2, e_3)$, we can consider a general basis of $\R^3$ and derive  directions therefrom.  
Secondly, in order to obtain a grid-free approach to safety analysis if the dimension is beyond 2, instead of constraining the GD on  planes, we can constrain the GD on more general regions, e.g. on the regions $x_j \in [i\delta, (i+1)\delta]$, $i\in \N, 0\leq i\leq \delta^{-1}, j\in\{1, 2, 3\}$.  \\
Note that once we select a plane (say $\phi$, or a region) it could happen that $\min_{x\in\phi}\mathbb{P}\left(\tau^{x}\geq T\right)<p$, which means that there is no intersection between $\phi$ and $\mathbb{A}^p_t$ and we must pass on to examining another plane (or region). 

\subsection*{Higher dimensions} 
We can apply the same reasoning on models with any dimension: namely, if we are in $\R^d$, then we can partition the considered region $A$ in sets of dimension $d-1$. 
Continuing this procedure we can go back recursively to the case $d=2$. 

\subsection*{Alternative approaches for higher dimensions}
An alternative grid-free approach is to ``explore the border'' without constraints that are relative to some sections, i.e. to generalize directly from the case $d=2$. 
Let us suppose that $x_\star$ is a point on the border; then we can compute $d-1$ orthonormal vectors $\{g^\star_1, \ldots, g^\star_{d-1}\}$ to $D_x \Prob(\tau^x\ge T)|_{x=x_\star}$, thus running the procedure recursively from any new point $x_\star \pm \gamma g^\star_j, j=1,\ldots,d-1$,  until we obtain a closed surface. 
However, attention is needed with the selection of the orthonormal points $\{g^\star_1, \ldots, g^\star_{d-1}\}$: indeed,  when $d>2$ there are infinitely many possibilities, but it would be convenient to find a possible ``orientation'' such that the exploration of the border is done in an orderly - see the discussion relative to Line~\ref{step:terminatewhile} of Algorithm~\ref{algo:walk_border_d2} above. 

\section{Experiments}\label{sec:exp}

In this section, we present a case study:  
the code can be found at \url{https://github.com/FraCose/Grid-free_prob_safety}. \\
For the experiments, we use a simulation-based approach, i.e. we use Monte Carlo (MC) techniques, and to reduce the variance we use antithetic Brownian paths \cite{Glasserman2003,Huynh2012}.  
\begin{remark}
We remark that the way $\mathbb{P}\left(\tau^{x}\geq T\right)$ and $D_x
\mathbb{P}\left(\tau^{x}\geq T\right)$ are computed it is not relevant for the idea presented in this work. Indeed, it is enough to be able to compute the quantities $\EE[\mathds{1}_{\{\tau\geq T\}}]$ and $\EE[\mathds{1}_{\{\tau\geq T\}}H_T]$ and plug them into the GD procedure. 
We refer to~\cite{Andersson2017, Frikha2019, HenryLabordere2017}  for the exposition of unbiased simulation methods. 
Other methods to compute these quantities are PDE techniques, which we expect to be computationally heavier.
\end{remark}
Before presenting the model for the case study, it is important to draw some general considerations on the discussed technique. 

\subsection*{Complexity} Let us recall the definition of $H_t$ and $\beta_t$: 
\begin{align}
H_T=&\int_{0}^{T}\frac{\mathds{1}_{\{t<\tau^{1}\}}}{\text{dist}(X_{t}, \partial A)^{2}}\beta_t\cdot dW_{t}, \\
\beta_t&=\sigma^{-1}(t,X_{t})\cdot J_{t}.
\end{align}
Computing $H_t$ can be expensive. 
To estimate the expectation via MC methods we use $N$ simulations and a time discretization step of $n^{-1}$, i.e. we split the time interval $[0,T]$ in $n$ steps.   
The stochastic processes to be simulated are %
$X_t, J_t, \beta_t, H_t$ and $\text{dist}(X_{t}, \partial A)$. 
The realization of the stochastic process $\beta$ has a total cost of $Nnd^3$, where $d^3$ is the cost related to the matrix inversion $\sigma^{-1}$, plus matrix multiplications. 
Moreover, an optimization problem to compute $\text{dist}(X_{t}, \partial A)$ has to be solved $Nn$ times. 
Nevertheless, we have to simulate $H_t, J_t$ only if $t<\tau_1$.
It is important to remark that we have analysed the computational cost of computing only one step of the gradient descent procedure, but many are necessary to converge and explore the space. 

If we are interested in a relatively low-dimensional problem the matrix inversion can be solved analytically, or leveraging special forms of $\sigma$, e.g. tri-diagonal, upper(lower)-triangular. 
This increases the stability of the procedure and reduces in part its complexity, although the overall complexity remains $Nnd^3$, being dominated by matrix multiplications. 
A second improvement is to consider particular forms for the region $A$ that can be advantageous for computing $\text{dist}(X_{t}, \partial A)$, e.g. a sphere, a parallelepiped or a simplex -- although non-smooth regions are not covered by the assumptions of this work. 
Furthermore, both the arguments just discussed allow the usage of GPU acceleration more easily, which ``artificially'' reduces the complexity in $N$.

\subsection*{Bias} It is important to remark that the steps done by the Gradient Descent algorithm are stochastic and biased. Indeed, we do not compute the exact probability $\Prob(\tau^x\geq T)$, but we discretize the time, therefore computing $\Prob(\tau_n^x\geq T)$; recall that in \cite{Gobet2000} it is shown that 
\[
|\Prob(\tau^x\geq T)-\Prob(\tau_n^x\geq T)|\leq O(n^{-1/2}),
\]
where $\tau_n$ represents the discrete stopping time of the Euler Scheme associated with Equation~\eqref{eq:SDE_ch6sec2}. 
Moreover, the algorithm is stochastic since we approximate $\Prob(\tau_n^x\geq T)$ using MC techniques. 
Therefore we have to consider that \cite{Gobet2000}
\begin{align}
\left|\mathbb{P}\left(\tau^{x}\geq T\right)\!-\!\widehat{\mathbb{P}\left(\tau_{n}^{x}\geq T\right)}\right|\leq&O(n^{-1/2})\!+\!O\left(\!\frac{1}{\sqrt{N}}\!\right)\!\!Z,
\end{align}
where the hat denotes an MC estimator of the quantity of interest, $Z$ represents a standard normal random variable, and $N$ is the number of simulations. 
A similar error bound might be derived for the other term $\left|D_{x}\mathbb{P}\left(\tau^{x}\geq T\right)\!-\!\widehat{D_{x}\mathbb{P}\left(\tau_{n}^{x}\geq T\right)}\right|$ \cite{Gobet2000,Gobet2007,Gobet2010}, though an adaptation is needed due to the presence of $\tau_1$ in the definition of $H$ in Theorem~\ref{th:Delta_Malliavin}. 
Due to these biases, we have noticed that reducing the variance helps the GD to converge better (cf. use of antithetic Brownian paths mentioned above): %
for instance, when the (norm of the) gradient $D_{x}\mathbb{P}\left(\tau^{x}\geq T\right)$ becomes small, the error could dominate and the gradient descent step could not work properly; this is especially the case when we simulate paths starting from points close to the border of $A$. 

\subsection*{Hyper-parameters} 
The hyper-parameters to be chosen for the procedure are the following: 
\begin{enumerate}[label=\textit{(\roman*)}]
\item $n$, time discretisation step - in principle the higher the better, but $n$ has a big impact on the computational time, since it cannot be parallelised. Through experiments, we have learnt to start with a relatively fine time discretisation step. 
\item $N$, Monte Carlo simulations - increasing $N$ reduces the variance of the MC methods. $N$ has a relatively low impact since the number of samples can be parallelised using a GPU.
\item $\lambda$, the ``learning rate'' of the GD procedure in  Equation~\eqref{eq:GD_step} - $\lambda$ must be chosen carefully. 
While we are doing the first minimization, i.e. while we are searching for a first point on $\partial \mathbb{A}^p_T$ (exploration phase), $\lambda$ can be quite high (more than $\num{1E-3}$, as suggested in \cite{Kingma2014}).  
Instead, if we are considering the minimization procedure in Algorithm~\ref{algo:walk_border_d2}, since we should be already close to the border we should select a small $\lambda$.
\item $\gamma$, the ``border exploration'' parameter in Algorithm~\ref{algo:walk_border_d2} -  
$\gamma$ indicates how fine-grained we wish the approximation of $\partial \mathbb{A}^p_T$ to be. 
If it is selected to be excessively small, the exploration of the border $\partial \mathbb{A}^p_T$ will be quite slow. 
\end{enumerate}

\subsection*{Acceleration of the exploration} Algorithm~\ref{algo:walk_border_d2} is a good starting point to explore the border $\partial \mathbb{A}^p_T$, however in practice care must be taken. For the following discussions, we consider the problem to be in a 2-dimensional space as a base case.  

\begin{figure}[tbh!]
    \centering
        \includegraphics[height=4.5cm,width=0.45\textwidth]{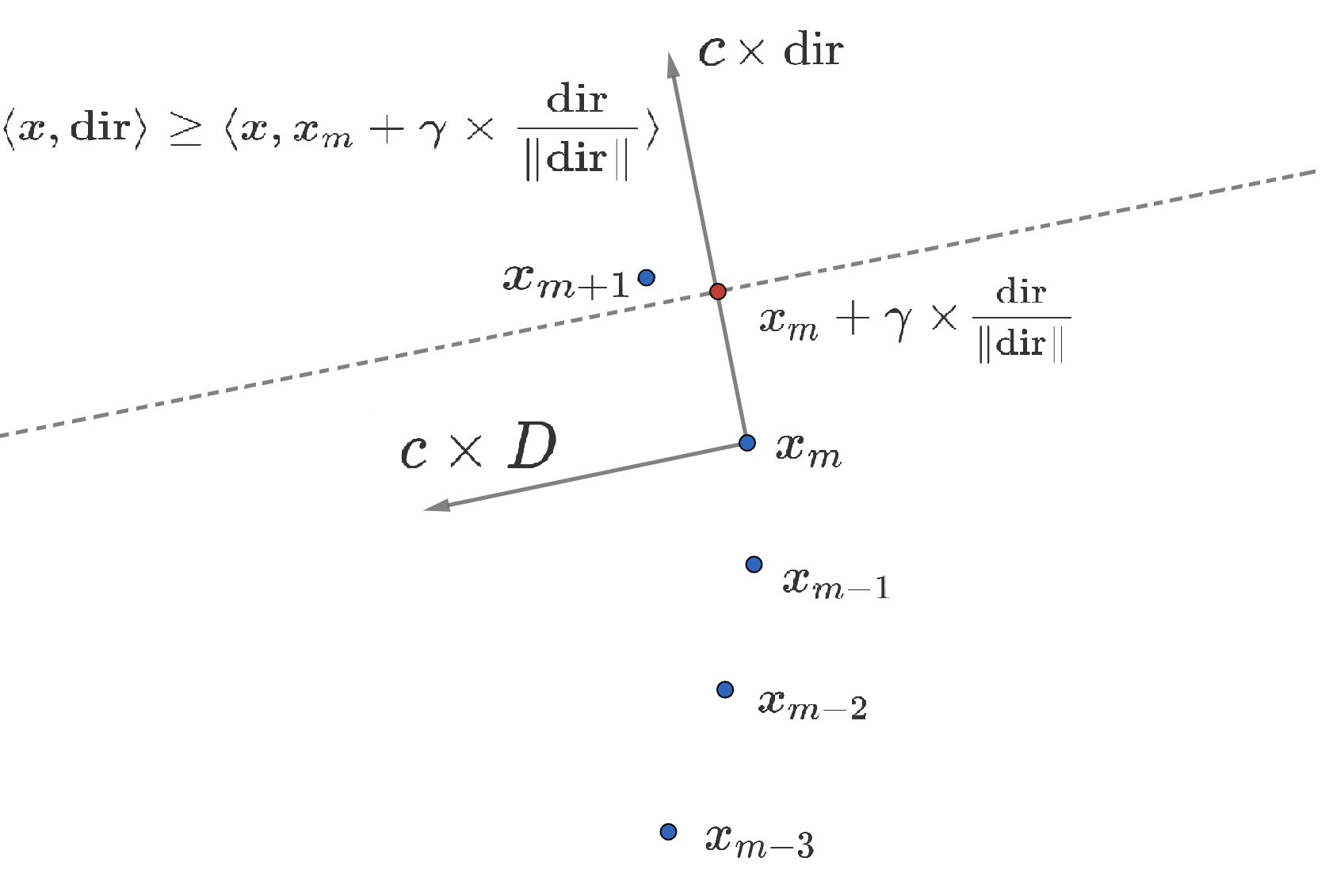}
    \caption{
    {Representation of how the algorithm explores the border: given the points $\{x_m, x_{m-1}, \ldots\}$ already found by the algorithm, it proposes the red point $x_m+\gamma\times \text{dir} / \|\text{dir}\|$ as the new point of the region, and from there it runs the GD to find the new point $x_{m+1}$. It is possible that on the half-plane where the algorithm looks for the new point does not exist a point of $\partial\mathbb{A}^p_T$, therefore it is necessary to update the constraint, see Algorithm~\ref{algo:fix_constraint}.}}
    \label{fig:directions}
\end{figure}

Firstly, we would like to explore with an orientation, e.g. clockwise, such that we do not go back to a region already explored. 
This can be done in principle, but sometimes the gradient approximation can be (quite) wrong, especially close to the border of the considered region $A$ or because the chosen discretization time step $n$ is too coarse. 
To solve this problem, we check if there are already points computed in the direction we are going to explore. However, selecting an ``optimal'' number of points is an open question that depends on the curvature of $\partial \mathbb{A}^p_T$, which a-priori is unknown. 
Another heuristic is to constrain the algorithm to search the new point on $\partial \mathbb{A}^p_T$ in a given region, see below and Figure~\ref{fig:directions}.
A more sophisticated alternative is to split the region into subspaces and search the border $\partial \mathbb{A}^p_T$ locally. This technique would also increase the level of parallelisation \cite{Suffern1990}.

Secondly, going towards the direction perpendicular to the gradient, see Section~\ref{sec:walk}, is only an approximation and sometimes, depending on the local curvature,  can be quite imprecise. 
To improve this approximation we have considered the following procedure. 
Let us imagine that we have computed a certain number of points on $\partial \mathbb{A}^p_T$, in order $\{x_1, x_2, ..., x_m\}$. We can compute the parabola equation (since the plane is fixed) that approximates the points $\{x_1, x_2, ..., x_m\}$, and use this equation as a second possible approximation. This can be thought an approximation of the second-order information of the curve $\Prob(\tau_n^x\geq T)-p=0$ in $x$. 
Later, we can choose the new direction as a weighted average of the perpendicular to $D_x \Prob(\tau_n^x\geq T)$ and the value of the approximated parabola $p(x)$ in $x = 2x_m-x_{m-1}$. 
As there are several ways to compute the weights, we use the past distances between the points found on $\partial\mathbb{A}^p_T$ and the forecasts relative to the gradient and the parabola approximation, see Figure~\ref{fig:directions} and the code for more insight. 
In this way, when the curvature of $\partial \mathbb{A}^p_T$ ``changes'' the algorithm starts following more closely the gradient (if the approximation error is low), otherwise it follows an average which experimentally is closer to the parabola forecast. 
Experimentally, this procedure accelerates the exploration, since it reduces the approximation error relative to the gradient.

Finally, it is better to constraint the space where the algorithm searches for the next point of the border. 
In Figure~\ref{fig:directions} it is shown how we proceed. 
Once one point $x_m$ on the border is found, i.e. $ \Prob(\tau_n^{x_m}\geq T) \approx p$, we compute the gradient ($D$) and the direction to follows\footnote{Possibly as a weighted average of the gradient and some local approximation of the curvature as explained before.} (dir). 
Given dir and $\gamma$, it is possible to search the new point only in the part of plane where there are not ``recent points'' considering the line perpendicular to the direction passing through the point guess $x_m+\gamma\times \text{dir}/\|\text{dir}\|$. \\
{It is possible that the constraint does not allow the optimization procedure to find a point $\Prob(\tau_n^{x_\star}\geq T) \approx p$, therefore }
if the solution of the GD returns, after a certain number of iterations, a point $x_\star$ s.t. $ \Prob(\tau_n^{x_\star}\geq T) \not\approx p$, 
then it is important to update the direction dir and the corresponding constraint. 
The candidate we have chosen for the updated direction is $2x_\star-x_m$, up to some re-scaling, but other choices are available. For example, 
we have experimented that selecting $2x_\star-x_m$ accelerates the procedure over the choice $x_\star-x_m$. 
Moreover, it is necessary to reduce the step exploration $\gamma$, such that we get closer to the point $x_m$ and by continuity of $\partial\mathbb{A}^p_T$ we will find the point sooner or later. 
In Algorithm~\ref{algo:fix_constraint} we present a pseudo-code of the procedure. 
\begin{algorithm}[tbh!]\caption{{Adaptive constraint for the GD procedure }}\label{algo:fix_constraint}
\begin{algorithmic}[1]
\State{Given $x_m$}
\State{Run the constrained GD from $x_m+\gamma\frac{dir}{\|dir\|}$ up to a point $x_\star$}%
\State{step $\leftarrow 1$}
\While{$\mathbb{P}\left(\tau^{x_\star}\geq T\right)\not\approx p$}
\State{dir $\leftarrow 2x_\star-x_m$}
\State{$\bar{x}\leftarrow x_{m}+\frac{\gamma}{2\cdot\text{step}}\cdot\frac{dir}{\|dir\|}$}
\State{Update the plane using the new dir and $\bar x$ }
\State{Run the constrained GD from $\bar x$ up to a point $x_\star$}
\State{step $\leftarrow \text{step} + 1$}
\EndWhile
\end{algorithmic}
\end{algorithm}
\begin{figure}[tbh!]
    \centering
        \includegraphics[height=4cm,width=0.31\textwidth]{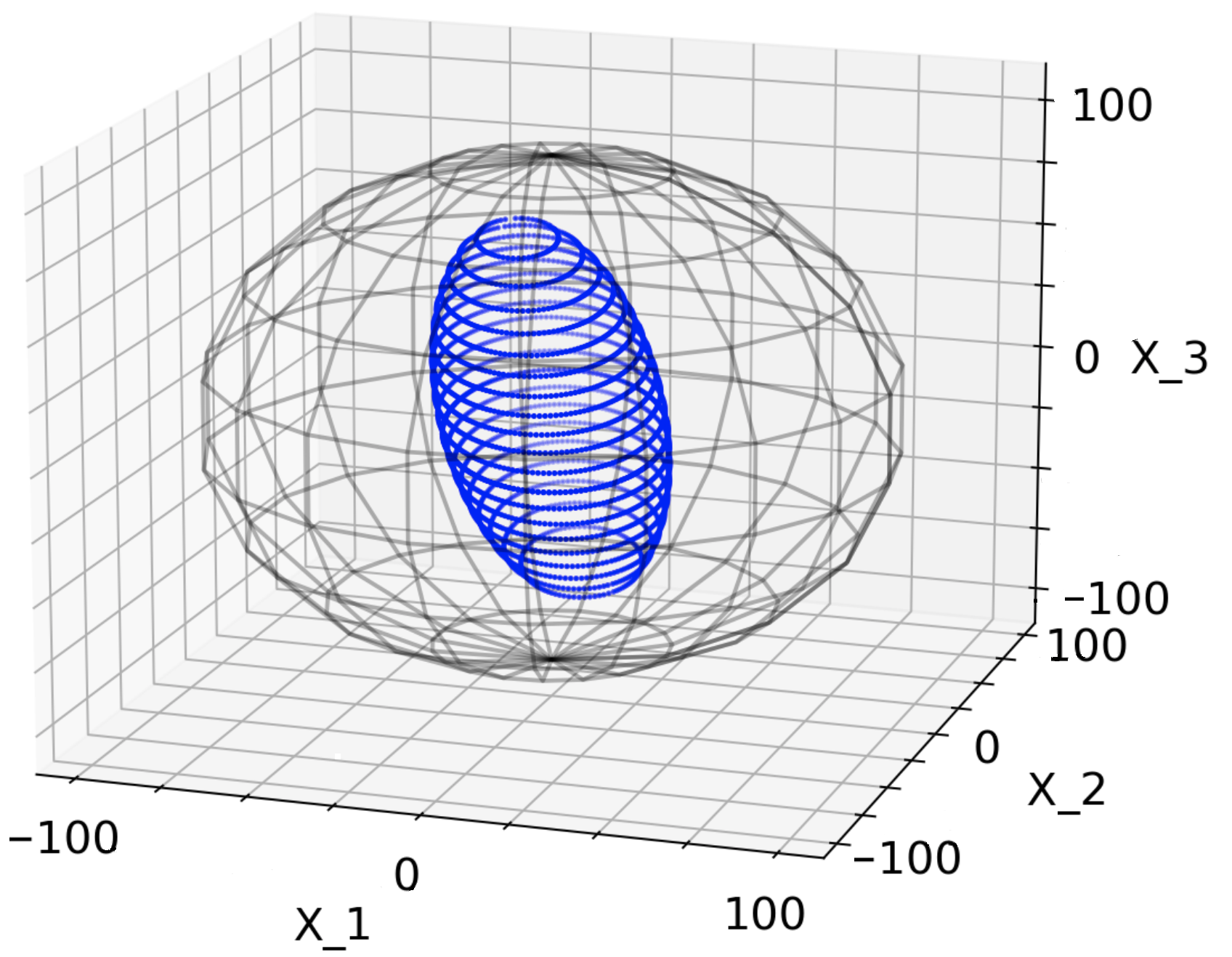}
        \includegraphics[height=4cm,width=0.31\textwidth]{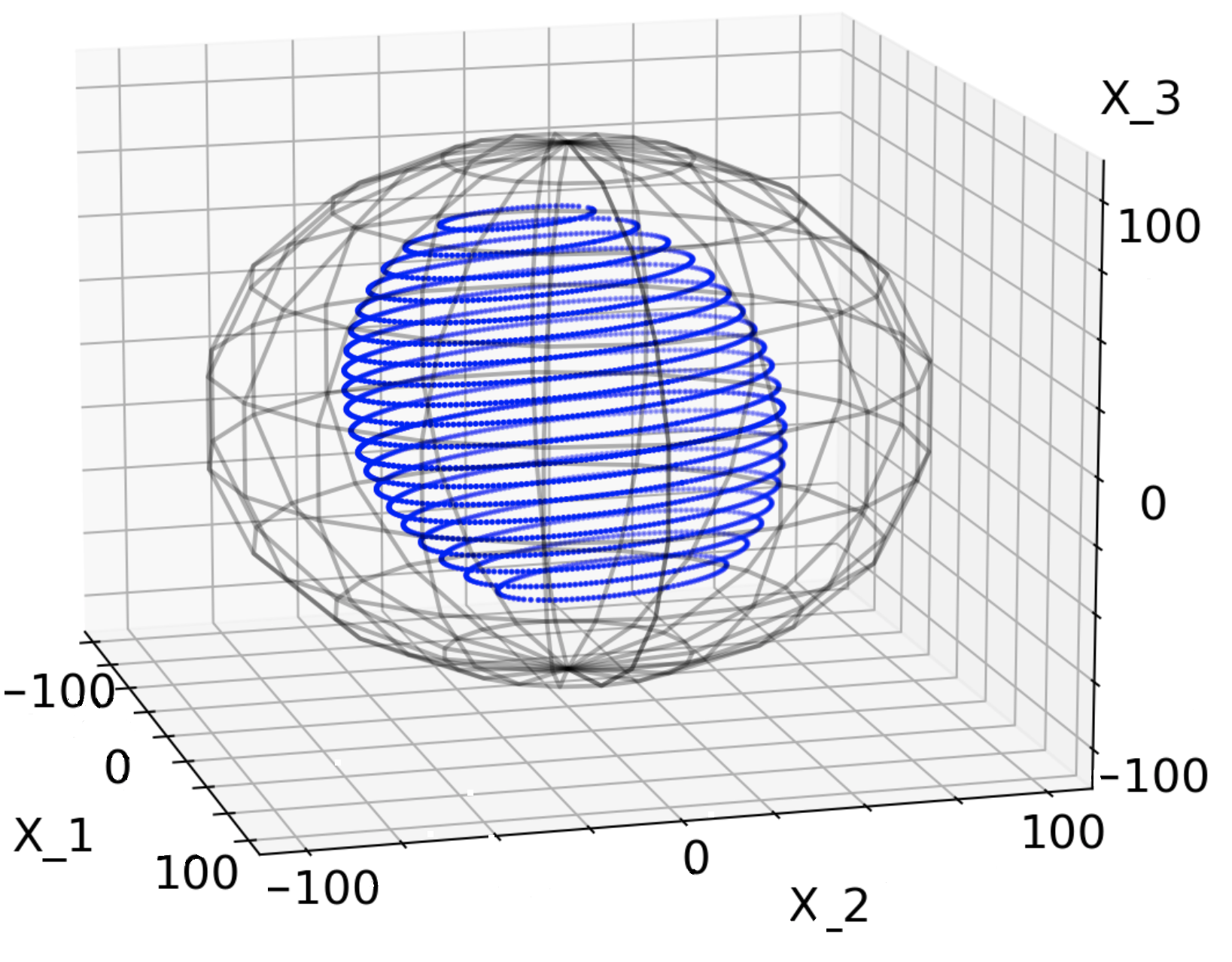}
        \includegraphics[height=4cm,width=0.31\textwidth]{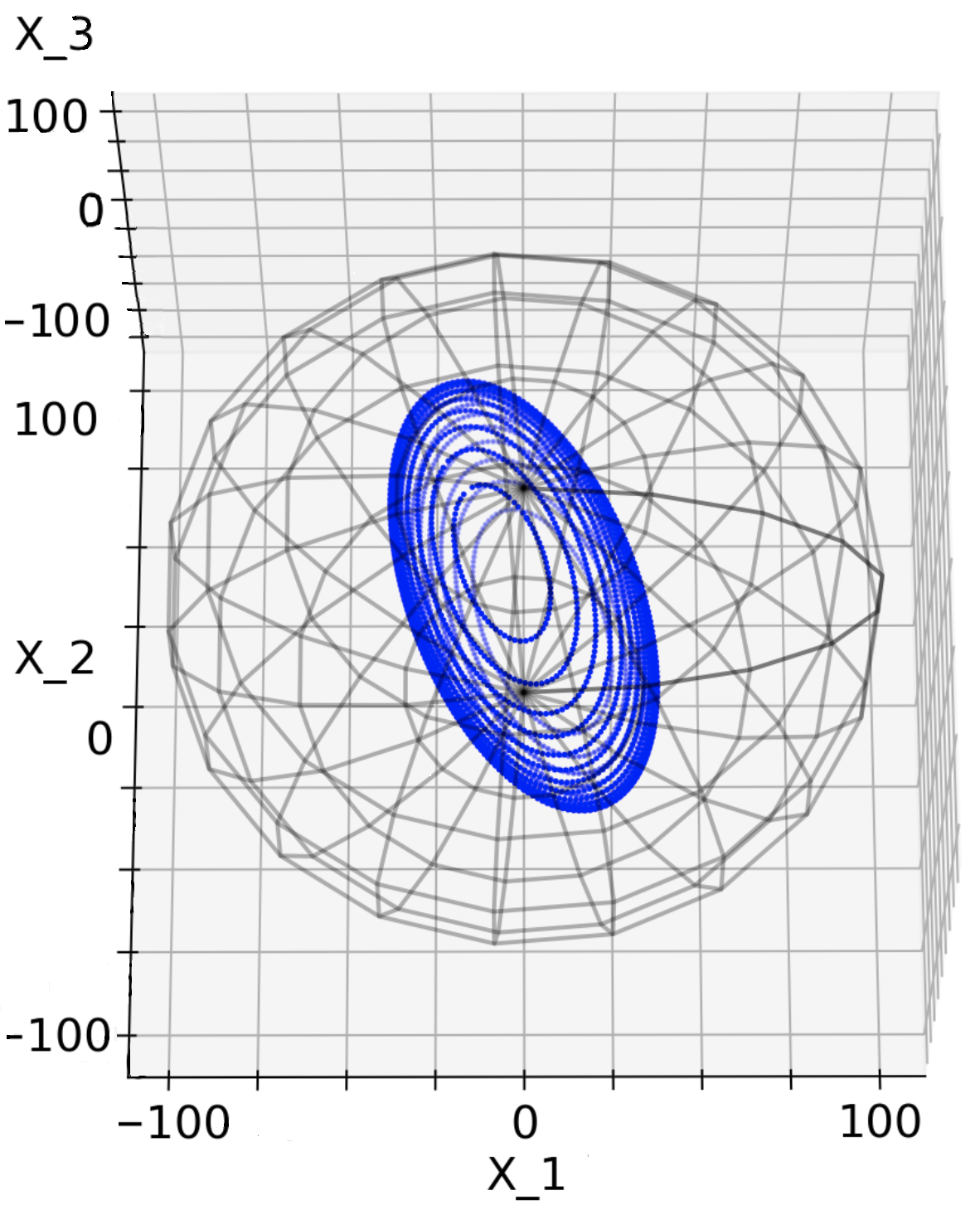}
    \caption{
    The plots show the surface $\partial\mathbb{A}^{0.5}_1$ (in blue) found when the region $A$ is a (black) sphere. 
    This has been computed sectioning the region $A$ across 2D planes (we have discussed at the end of Section~\ref{sec:walk} the use of alternative grid-free approaches for exploration in 3D (and higher-dimensional) cases).
    }
    \label{fig:experiments_full}
\end{figure}
If the direction guides towards points already explored recently, because for instance the discretisation error is too high or due to the constrained updates in Algorithm~\ref{algo:fix_constraint},  given that on a plane the perpendicular vectors to a vector are two, it is enough to invert the direction. 

\subsection*{Case study} 

The model considered for the experiment is:
\begin{align}\label{eq:toy_experiment}
\!\!\!\!\!\!\!\left(\!\!\!\begin{array}{c}
dX_{t}^{(1)}\\
dX_{t}^{(2)}\\
dX_{t}^{(3)}
\end{array}\!\!\!\right)\!=\!&\left(\!\!\!\begin{array}{c}
X_{t}^{(1)}\\
\dfrac{1}{2}X_{t}^{(1)}+\dfrac{1}{2}X_{t}^{(2)}\\
\dfrac{1}{3}X_{t}^{(1)}+\dfrac{1}{3}X_{t}^{(2)}+\dfrac{1}{3}X_{t}^{(3)}
\end{array}\!\!\!\right)dt +\!\frac{1}{3}\left(\!\!\!\begin{array}{ccc}
\omega_{1},\quad & \omega_{2},\quad & \omega_{2}\\
\omega_{2},\quad & \omega_{1},\quad & \omega_{2}\\
\omega_{2},\quad & \omega_{2},\quad & \omega_{1}
\end{array}\!\!\!\right)dW_{t},\\
\omega_{1}:=&2\sqrt{1-\rho}+\sqrt{1+2\rho},\\\omega_{2}:=&-\sqrt{1-\rho}+\sqrt{1+2\rho}.
\end{align}

If we define $d\tilde{W}_t = \sigma dW_t$, where $\sigma$ is the diffusion matrix in Equation~\eqref{eq:toy_experiment}, then we have that $\operatorname{Corr}(d\tilde{W}^{(i)}_t, d\tilde{W}^{(j)}_t)=\rho$, $i\not= j$ and $i, j \in \{1, 2, 3\}$. In the experiment we have used $\rho=0.5$. 
For the region $A$, we have considered two cases:
a sphere with center at the origin and radius equal to $100$
and a cube with vertices between $(-100,-100,-100)$ and
$(100, 100, 100)$.
Note that in the second experiment ($A$ being a cube) the assumptions of the theoretical part of this work are not satisfied. Nevertheless, the procedure is still able to explore the border. 
\begin{figure}[H]
    \centering
        \includegraphics[height=5.5cm,width=0.47\textwidth]{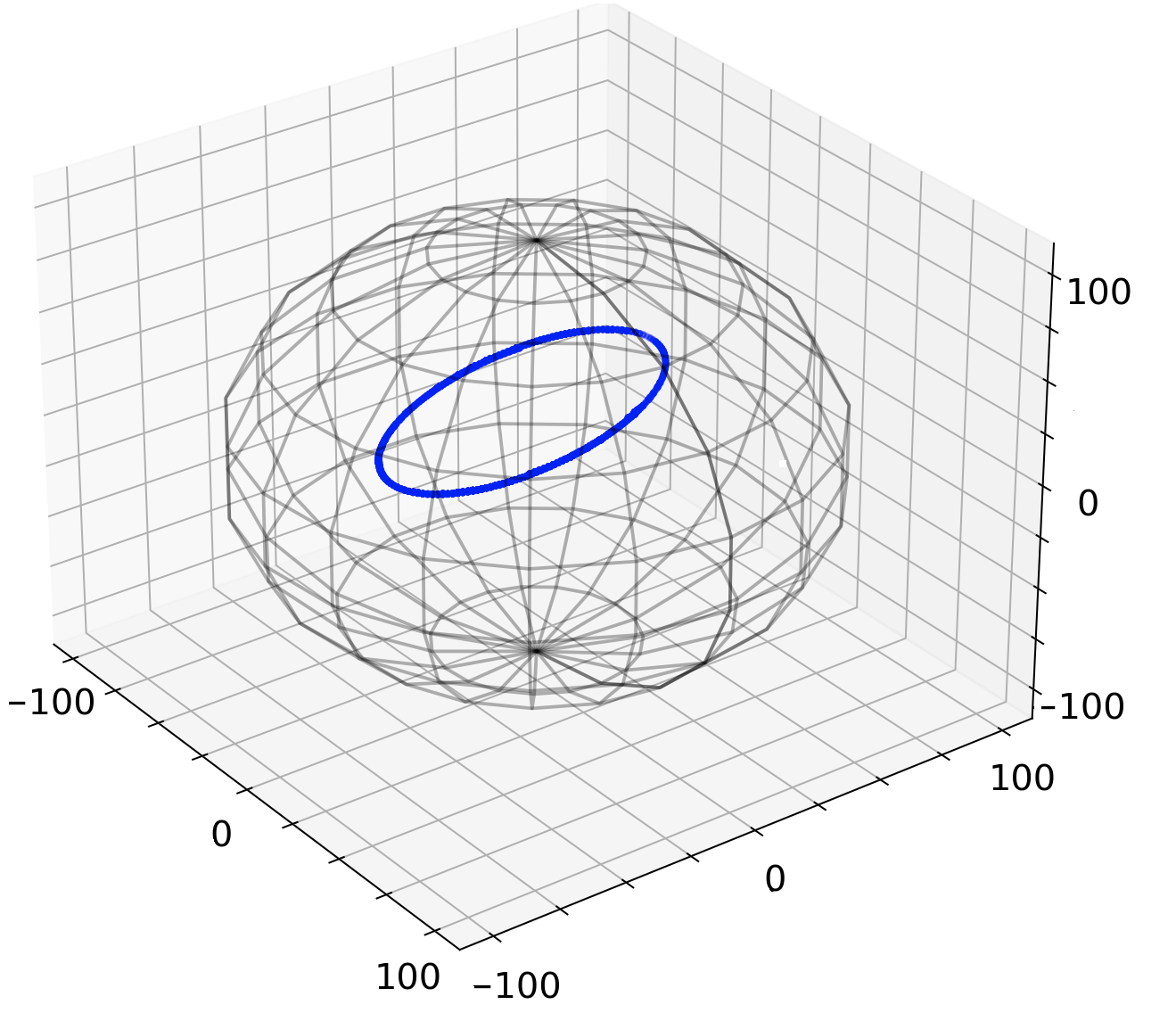}
        \includegraphics[height=5.5cm,width=0.47\textwidth]{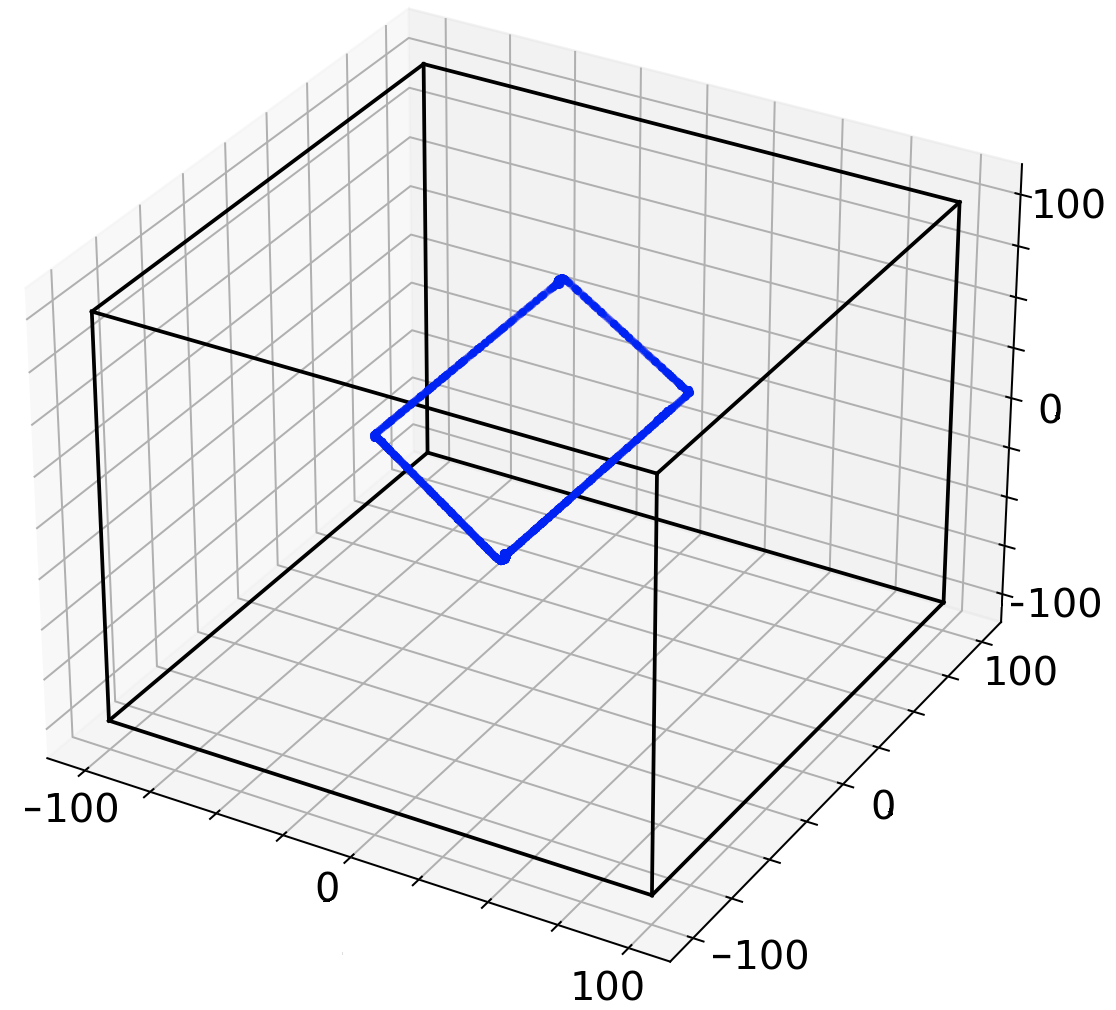}
        \includegraphics[height=5.5cm,width=0.47\textwidth]{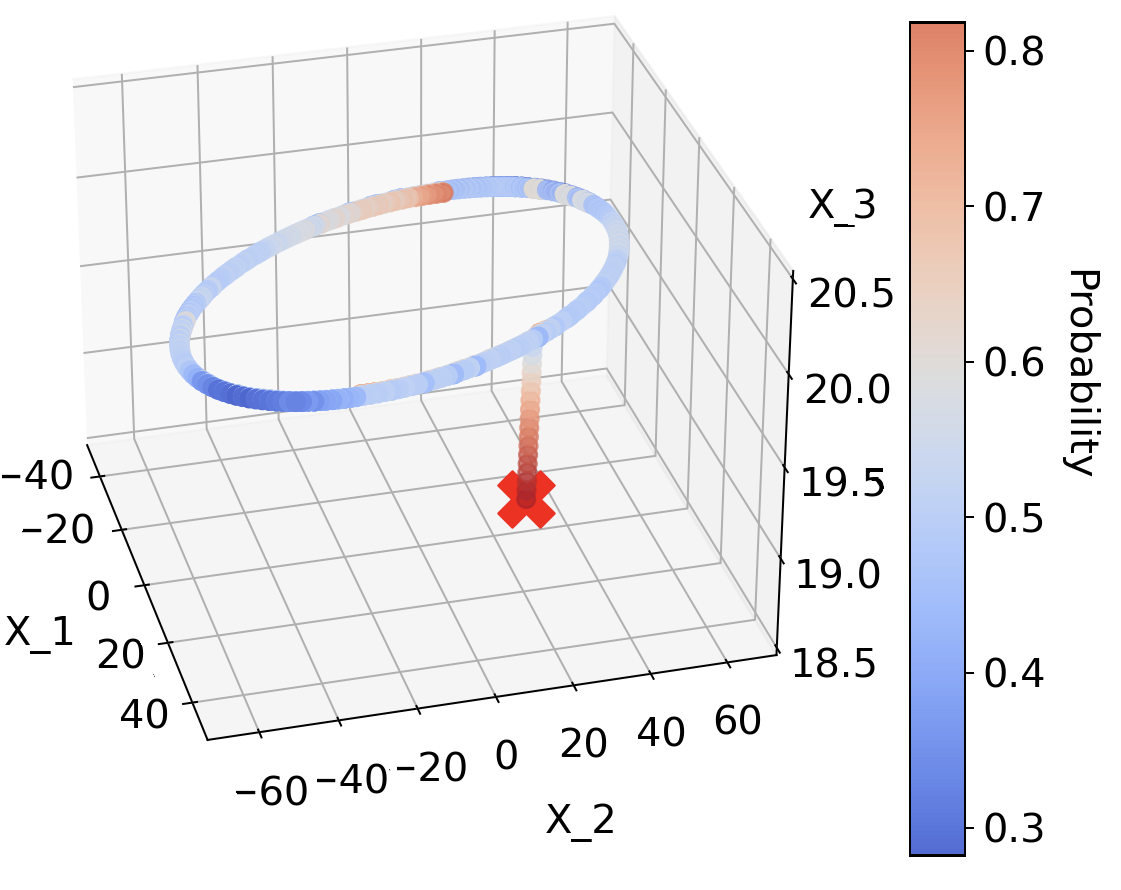}
        \includegraphics[height=5.5cm,width=0.47\textwidth]{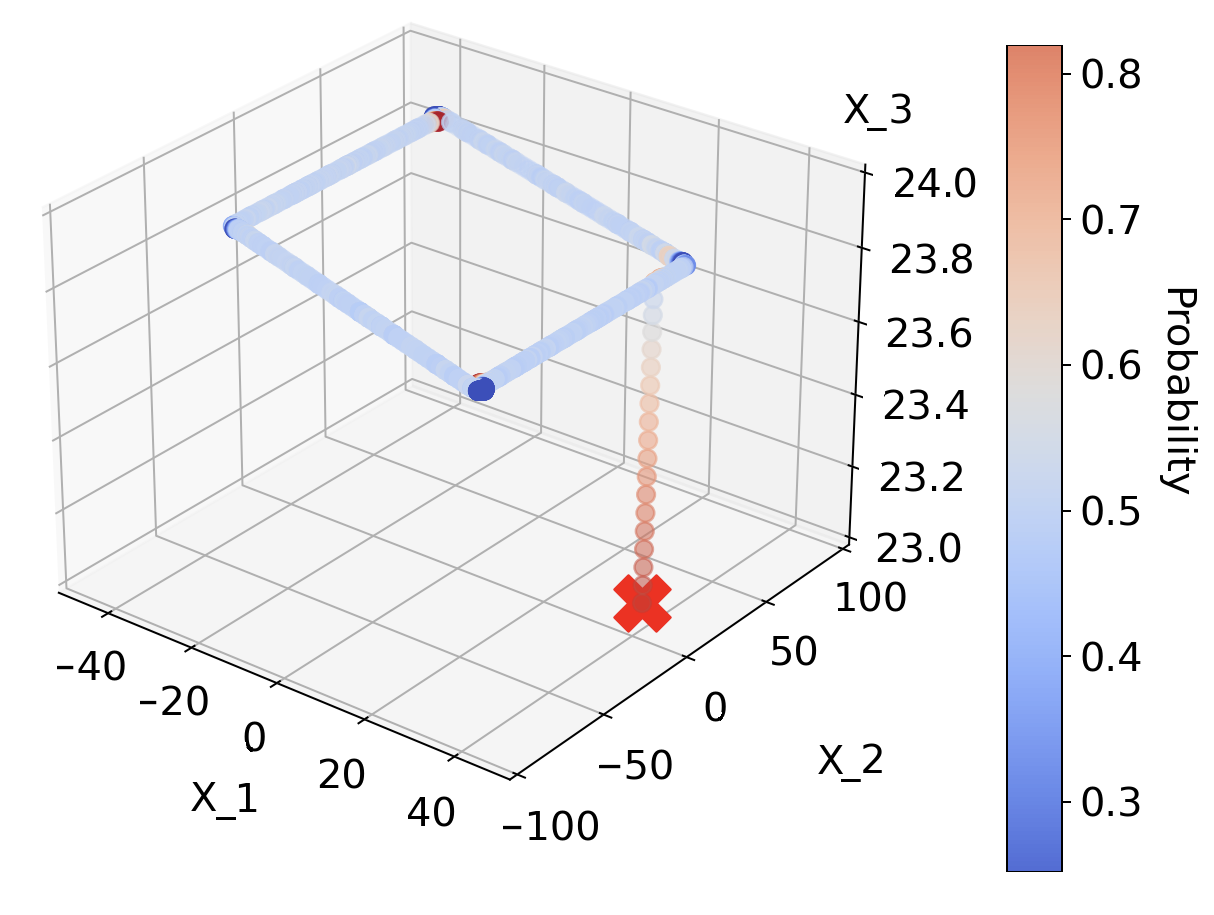}
        \includegraphics[height=5.5cm,width=0.47\textwidth]{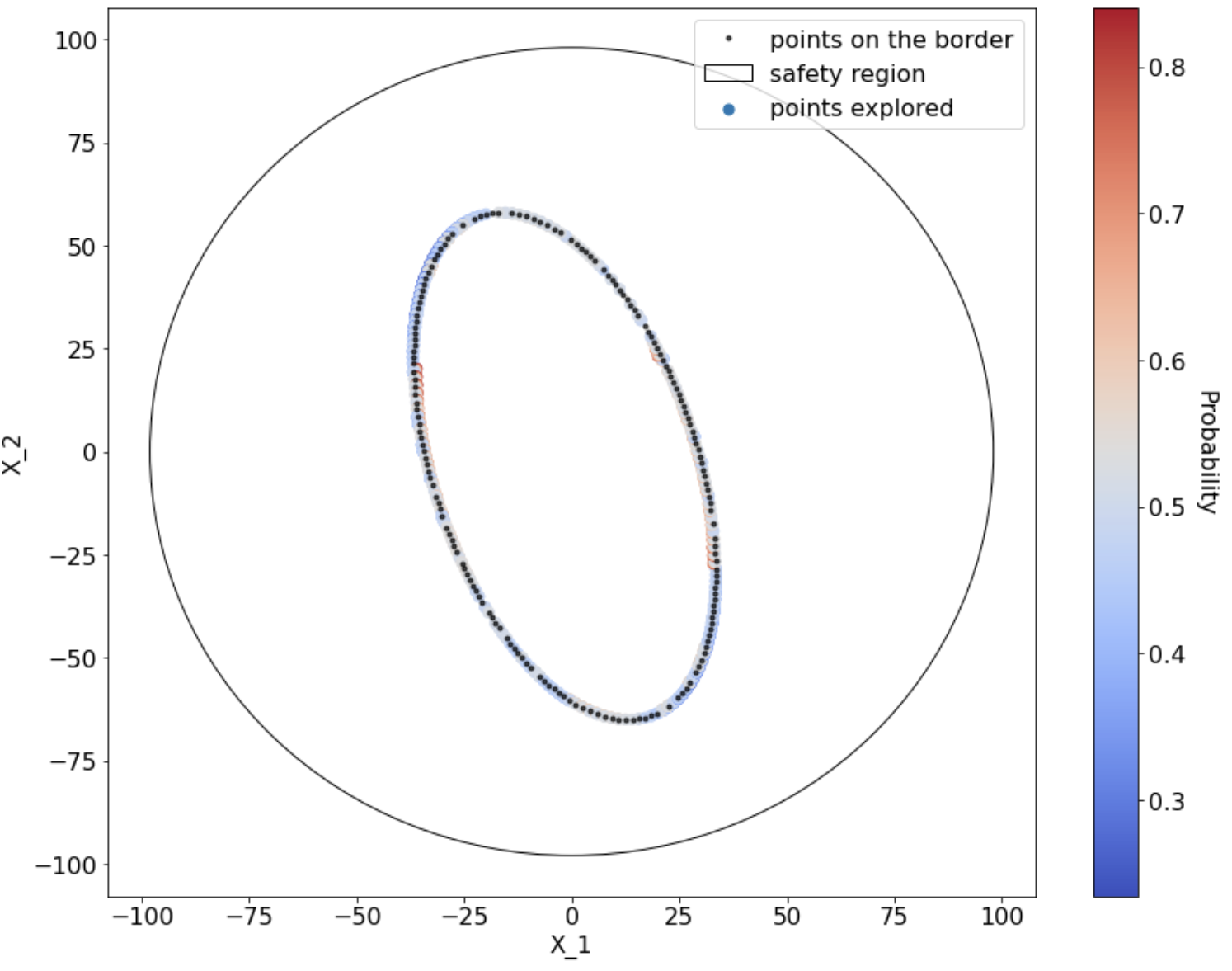}
        \includegraphics[height=5.5cm,width=0.47\textwidth]{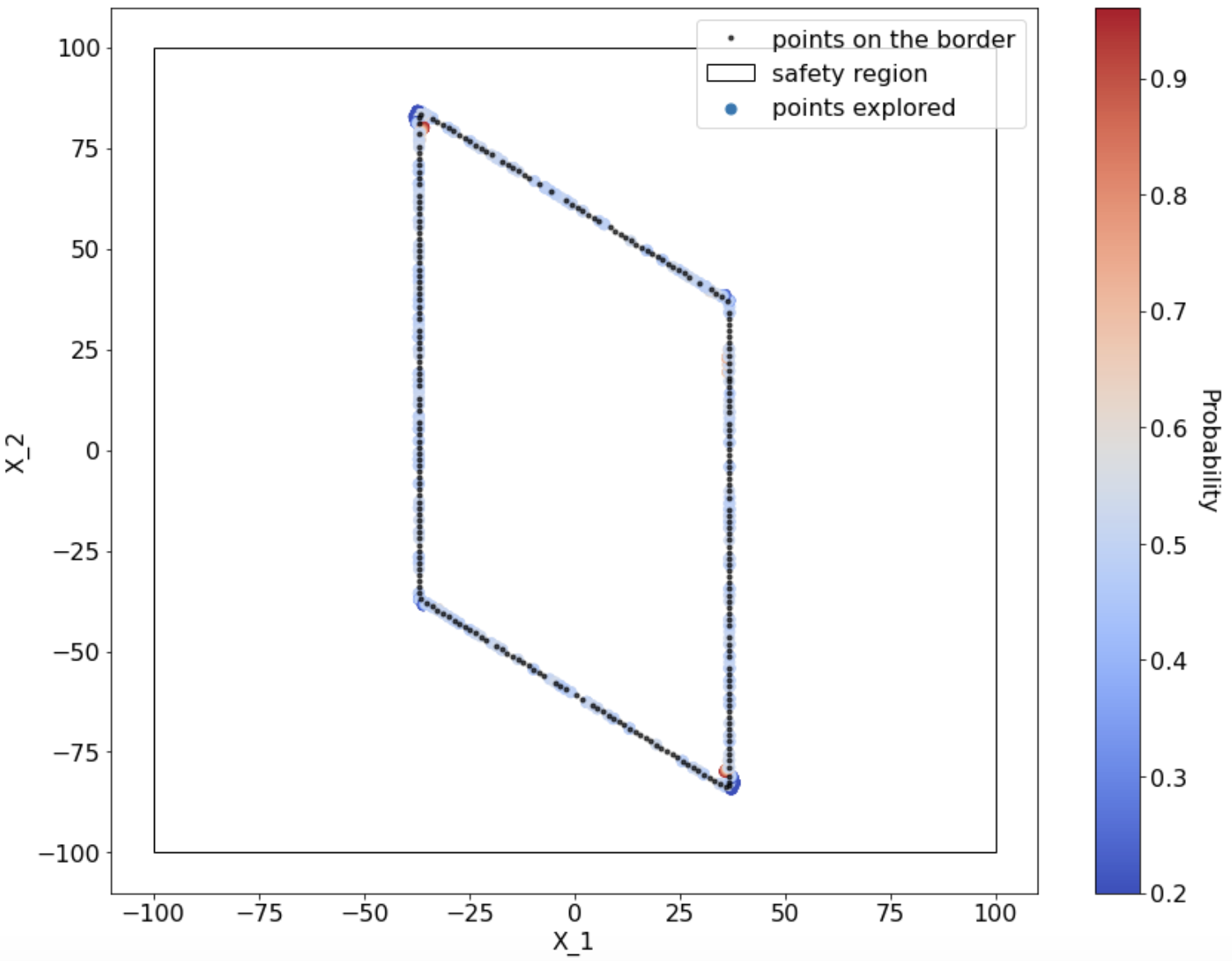}
    \caption{Left. The plots show the points found while and exploring the border of the 3D region $\mathbb{A}^{0.5}_1$, with respect to one plane (or section), when $A$ is a sphere.
    Top: full-view of the problem. 
    Middle: same plot but zoomed, with colour scale  showing the safety probabilities of the points explored during exploration. 
    It can be seen how starting from the red-cross point, with a probability approximately $0.8$, we arrive at the points with the desired probability $0.5$. %
    The presence of points with different colours to the one corresponding to the desired probability $0.5$ means that in those regions the GD has explored adjacent points. 
    Bottom: planar section of the space considered, fixing direction $x_3$. The black points represent those we consider being part of the border (up to an approximation error), whilst the circle points are those explored via the GD method. It can be seen that we explore points associated with probabilities between $0.2$ and $0.8$.\\
    Right. The plots show the points found seeking and exploring the border of the 3D region $\mathbb{A}^{0.5}_1$, with respect to one plane (or section), when $A$ is a 3D cube.
    }
    \label{fig:experiments}
\end{figure}
We consider the problem of computing 
$
\mathbb{A}^{0.5}_1 = \{x\in A : \Prob(\tau^x\geq 1)\geq 0.5\}.
$
We start at a point $x_0$ where $\Prob(\tau^{x_0}\geq 1)\not= 0.5$, then we minimize 
$\frac{1}{2}(\Prob(\tau^x\geq 1)- 0.5)^2$ until we obtain a point $x_\star$ s.t. $|\Prob(\tau^{x_\star}\geq 1)-0.5|<\text{err}$ -- in this case err represents the approximation errors due to the computation of $\Prob(\tau^{x_\star}\geq 1)$. 
From $x_\star$ we fix $x_\star^{(3)}$ and start Algorithm~\ref{algo:walk_border_d2}, i.e. we fix the plane $x^{(3)} = x_\star^{(3)}$, see Figure~\ref{fig:experiments} for the results of the experiments.
See Figure~\ref{fig:experiments_full} for a 3-dimensional representation of (possibly a portion of) $\partial\mathbb{A}^{0.5}_1$ in the case $A$ is a sphere. 
Instead of using plain Gradient Descent, we use ADAM \cite{Kingma2014}, a version of GD with momentum and adaptive learning rate that has shown impressive performance in Machine Learning and it is considered the state-of-the-art optimisation tool.  In particular, we prefer to include momentum, because we do not know whether $\frac{1}{2}(\Prob(\tau^x\geq 1)-0.5)^2$ is convex as a function of $x$.

The hyper-parameters chosen are $N =$ \num{10000}, $n = 200$, \textit{maximum iteration of the GD} (any time we use it) $= 50$, $\lambda =$ \num{5E-2}, $\gamma = 1.5$. With reference to the previous discussion on the approximation of the second order information of $\partial\mathbb{A}^p_T$, in order to compute the new direction, i.e. dir in Figure~\ref{fig:directions}, in addition to the gradient information, we use also the parabola approximating the previous $4$ points found on the border. For more information, we refer the reader to the code at \url{https://github.com/FraCose/Grid-free_prob_safety}. 

\section{Conclusions}
We have presented a new approach to find and compute probabilistic safety regions for stochastic differential equations (SDE) without resorting to the discretisation of their state space, which is by and large the standard approach in literature, which comes with limits related to precision and computational scalability. 
This is done by formulating an optimisation problem: to solve this, we have borrowed techniques and ideas from Malliavin Calculus and Mathematical Finance. 
We have discussed two formal results that allow one to explore relevant parts of the regions of interest, thus focusing computational load related to probabilistic safety computation for continuous-space models, such as SDEs. 
We have discussed possible algorithmic issues related the procedure, and offered strategies to cope with them.  
We conclude suggesting that more work on the generalisation of the approach on high-dimensional models in a completely automatic fashion is a goal worth pursuing. 

\subsection*{Acknowledgements and Disclosure of Funding}
The authors want to thank The Alan Turing Institute and the University of Oxford for the financial support given. 
FC is supported by the University of Oxford and The Alan Turing Institute, TU/C/000021, under the EPSRC Grant No. EP/N510129/1. 
HO is supported by the EPSRC grant ``Datasig'' [EP/S026347/1], The Alan Turing Institute, the Oxford-Man Institute and the University of Oxford.

\bibliographystyle{plain}
\bibliography{Biblio_greeks}

\end{document}